\documentclass{article}
\usepackage{graphicx,type1cm,eso-pic,color}
\usepackage{amsfonts, amssymb}
\newcommand{\dsp}{\displaystyle}

\newcommand{\wt}{\widetilde}

\newtheorem{theorem}{Theorem}[section]
\newtheorem{lemma}[theorem]{Lemma} 
 
\newtheorem{corollary}[theorem]{Corollary} 
\newtheorem{remark}[theorem]{Remark}

\title{\bf Parametric Resonance in Wave Maps}
\author{Tatsuo Nishitani and Karen Yagdjian\\
{}\\
\small Department of Mathematics, Osaka University, \\
\small Machikaneyama 1-1, Toyonaka, 560-
\small 0043, Osaka, Japan,\\ \small  nishitani@math.sci.osaka-u.ac.jp \\ 
{}\\
\small Department of Mathematics, 
\small University of Texas-Pan American, \\
\small  1201 W.~University Drive, 
Edinburg, TX 78541-2999,  
USA,\\  \small 
 yagdjian@utpa.edu }

\begin{document}

\date{}
\maketitle
\thispagestyle{empty} 
\vspace{-0.3cm}

\medskip

\setcounter{equation}{0}
\pagenumbering{arabic}
\setcounter{page}{1}
\thispagestyle{empty}

\begin{abstract}
In this note we concern with the wave maps from the  Lorentzian   manifold with the  periodic in time metric into the Riemannian  manifold,
which belongs to the 
one-parameter family of Riemannian  manifolds. That family   
contains as a special case  the   Poincar{\'e}  upper half-plane model. 
Our interest to such maps is motivated with  some particular type of the  Robertson-Walker spacetime
arising in the cosmology. We show that  small  periodic in time perturbation of the Minkowski metric generates  parametric resonance phenomenon. 
We prove that,
the global in time solvability in the neighborhood of 
constant solutions is not a stable property of the wave maps. \\
{\it Key words: wave maps, Floquet theory, parametric resonance, global solutions} 
\end{abstract}


\section{Introduction}
\label{S1}

In this note we prove instability,  under periodic in time perturbation of the Minkowski metric, of 
the global in time solvability  of the wave map equation in the neighborhood of 
stationary solutions.

Let $(V,g)$ be Lorentzian   and $(M,h)$ be  Riemannian manifolds   of dimensions $n+1$ and $d$, respectively. 
Let $u$ be a continuous mapping from
$V$ into $M$:
\[
u\,:\,\, V \longrightarrow M.
\]
The wave map is  the mapping  $u$ that satisfy the Euler-Lagrange equations for the Lagrangian   
\[
{\mathcal L}(u):= \int_{V} \frac{1}{2} \langle d u, d u \rangle_{T^*V\otimes u^{-1}TM} dV
\]
with $dV = \sqrt{g }dx_1 \wedge  \ldots \wedge  dx_{n+1} $ in local coordinates, being the volume form of V. 
From now on the Einstein's summation
convention is in force. Written in the local coordinates on $V$ and $M$ the  Euler-Lagrange  equations read 
\begin{eqnarray}
\label{WME}
&  & 
g ^{\alpha \beta } \left(  \partial_{\alpha\beta}^2  u^A - \Gamma _{\alpha  \beta }^\lambda  \partial_\lambda   u^A 
 + \Gamma _{BC }^A  \partial_\alpha    u^B  \partial_\beta     u^C \right)=0\,.
\end{eqnarray}
Here $  \Gamma _{\alpha  \beta }^\lambda$ and $\Gamma _{BC }^A $ denote  the components of the Riemannian 
connections of $g $ and $h $, respectively. The wave map equation is invariant under isometries of  $(V,g)$ and $(M,h)$.
\smallskip

We are interested in the Cauchy problem associated with (\ref{WME}). More precisely, for given initial data 
$(u(0),\partial_t u(0)) \,:\, S_0 \longrightarrow M\times TM$ 
at time $t=0$, we look for a wave map $u$ extending these globally in time. First, consider the case of  $(V,g)$ being 
the standard Minkowski space equipped with metric
\begin{eqnarray}
\label{MinMetr}
&  &
g ^{\alpha \beta } =\mbox{\rm diag}(-1,1,1,\ldots,1) \,.
\end{eqnarray}
In connection with this metric we recall conjecture of Klainerman and related known results for  the   Poincar{\'e}  upper half-plane 
model (the standard hyperbolic plane).
\smallskip

\noindent
{\bf Conjecture} (Klainerman \cite{Krieger_2004}) {\it Let $({\mathbb H}^2,h)$ be the standard hyperbolic plane. Then classical wave maps originating on ${\mathbb R}^{2+1} $ exist for  arbitrary smooth initial data.}
\smallskip

The following partial result towards the conjecture has been established in \cite{Krieger_2004}.

\begin{theorem} \mbox{\rm (Kreiger~\cite{Krieger_2004})} \,\,Let  $({\mathbb H}^2,h)$, 
${\mathbb H}^2:=\{(u^1,u^2) \in {\mathbb R}^2 \,|\, u^2>0 \}$, be  the standard hyperbolic plane 
with the metric tensor  $\dsp h_{ij}du^i du^j = \frac{1}{(u^2)^2}\left( (du^1)^2+(du^2)^2\right)$.
Then given initial data $u[0] : \{0\}\times  {\mathbb R}^2 \longrightarrow {\mathbb H}^2\times T{\mathbb H}^2$ 
which are sufficiently small in the sense that 
\begin{eqnarray}
\label{epsilon}
\int_{\{0\}\times  {\mathbb R}^2}\sum_{\alpha =0}^2 \left( \left[ \frac{\partial_\alpha u^1}{u^2} \right]^2+  \left[ \frac{\partial_\alpha u^2}{u^2} \right]^2\right) dx < \varepsilon 
\end{eqnarray}
for suitably small $\varepsilon >0$, there exists a classical wave map from ${\mathbb R}^{2+1} $  to $ {\mathbb H}^2 $ extending
these globally in time.
\end{theorem}
In particular, the global wave map exists for any small, in the sense of (\ref{epsilon}),  
compactly supported perturbations of the initial data of the constant wave map, which has
vanishing integral (\ref{epsilon}). 
In other words, {\it the global solvability of the wave map equation is a stable property in the small neighborhood of the global constant solutions.}
\smallskip

The   answer to the Klainerman's conjecture  as well as the scattering result for the wave map  are given in \cite{Krieger-Schlag}.
In particular, it is proved in \cite{Krieger-Schlag}
that if $M$ is a hyperbolic Riemann surface, and initial data $ (u(0),\partial_t u(0)) \,:\, S_0 \longrightarrow M\times TM$ 
are smooth and $u(0)=const$, $ \partial_t u(0)  =0$  outside of some compact set, then the wave map evolution $u $ 
of these data as a map ${\mathbb R}^{2+1} \longrightarrow M$ exists globally as a smooth function.
\smallskip

In this paper  we are interested in the case of the  Riemannian manifold   $(M,h)$  which  belongs to one-parameter family of manifolds containing  
the Euclidean half-space and  the Poincar{\'e} upper half-plane model  $({\mathbb H}^2,h)$.
In fact, that family consists of the Riemannian manifolds, which are 
 the half-plane $ \{(u^1,u^2) \in {\mathbb R}^2 \,|\, u^2>0 \}$ equipped with the metric   $\dsp h_{ij}du^i du^j = \frac{1}{(u^2)^l}\left( (du^1)^2+(d u^2)^2\right)$,
where the parameter $l$ is a real number. For $l=0$  the metric is Euclidean, while for $l=2$ it is the metric of the standard hyperbolic plane.
Those are the only two  manifolds of this family  which have constant curvature. 
\smallskip

In the present paper we examine the stability of the global solvability of the wave map equation with respect to the perturbation of the metric $g$.
First, we  prove that the only stationary solutions of the equation (\ref{WME}) are the constant solutions.  
Then, we  show that  
the global in time solvability can be destroyed     
by parametric resonance phenomena. (For the scalar quasilinear wave equation  
it was proved in \cite{KY_Torino},\cite{YagJMAA2001}.) 
Thus, {\it the small data global solvability is unstable with respect to the arbitrary small periodic perturbation of the metric tensor $g$.}
More precisely, we prove that the local solution obtained by small, for given $G \in {\mathbb N} $  in the sense of the following integral 
\begin{eqnarray}
\label{epsilonG}
\int_{\{0\}\times  {\mathbb R}^n}\sum_{|\gamma | =1}^G  \left( \left[ \frac{\partial^\gamma  u^1}{(u^2)^{l/2}} \right]^2+  \left[ \frac{\partial^\gamma u^2}{(u^2)^{l/2}} \right]^2\right) dx  ,
\end{eqnarray}
 smooth compactly supported perturbations
of the initial data of the global solution, in general, cannot be extended 
 globally in time. For the parametric resonance phenomena in the scalar wave map-type  hyperbolic equations see \cite{yagdjian_birk} and references therein. 
Then, according to  \cite{Ueda} (see  also   references therein)   the parametric resonance phenomena in the linear scalar wave equations can be localized in the space.     
  \smallskip

The Cauchy problem for  the wave maps in the perturbed Minkowski spacetime is considered in \cite{Choquet-Bruhat_ND}.
More precisely, assume that $V= S\times {\mathbb R}$, with $S$ an $n$-dimensional orientable smooth manifold, and let $g$ be a Robertson-Walker metric  
$
g   = -dt^2+R^2(t) \sigma$,  
where $\sigma =\sigma _{ij}\,dx^i \,dx^j $ is given smooth time independent metric on $S$, with non-zero injectivity radius. 
The Christoffel symbols of the Robertson-Walker metric are
\begin{eqnarray*}
&  &
\Gamma ^i_{jh} =\gamma ^i_{jh} , \quad \Gamma ^0_{00} = \Gamma ^0_{j0} =\Gamma ^i_{00} =0 , \quad 
\Gamma ^i_{0j} =R^{-1}R'\delta ^i_j , \quad \Gamma ^0_{ij} = RR'\sigma _{ij} , \quad R':=\frac{d R}{d t},
\end{eqnarray*}
where $\gamma  $ denotes the Christoffel symbols in the metric $ \sigma $.
The following result is known.

\begin{theorem}\mbox{\rm (Choquet-Bruhat~\cite{Choquet-Bruhat_ND})}
Let $(S\times {\mathbb R},g)$ be a Robertson-Walker expanding universe with the metric 
\begin{eqnarray}
\label{R-W_Metr}
g   = -dt^2+R^2(t) \sigma, 
\end{eqnarray}
with $(S,\sigma ) $ a smooth Reimanian manifold with  non-zero injectivity radius, of dimension $n \leq 3$, and $R$ a positive increasing function of $t$ such that 
$1/R(t)$ is integrable on $[t_0,\infty)$.

Let $(M,h)$ be a proper Riemannian manifold regularly embedded in ${\mathbb R}^N$ such that \, Riem($h$) is uniformly bounded. 

Then there exists a global wave map from $(S\times [t_0,\infty),g ) $ into  $(M,h)$ taking Cauchy data $\varphi  $, $ \psi $ with $D\varphi  $ and $\psi  $ in $H_1$
 if the integral of $1/R(t)$ on $[t_0,\infty)$ is less than some corresponding number $M(a,b)$.
\end{theorem}
The number $M(a,b) $ depends on the initial data. Thus, (see Corollary on page 45~\cite{Choquet-Bruhat_ND}) under hypothesis of the theorem, for any finite value of the integral of $1/R(t)$ on $[t_0,\infty)$ 
there is an open set $U$ of   initial data in $H_1\times H_1$  such that if $(D \varphi , \psi ) \in U$, then there exists a global wave map taking the Cauchy data $(\varphi , \psi )$.
In particular, this is true for the de~Sitter model of universe with  $R(t)=\exp(\Lambda t)$, $\Lambda >0 $.
  \smallskip

In Robertson-Walker geometry of positive curvature the space time interval is  (see Sec.9.8~\cite{Ohanian-Ruffini}, and p.131~\cite{Hawking})
\[
ds^2= -dt^2+a^2(t)\left(  d\chi ^2+\sin^2\chi (d \theta ^2 +  \sin^2 \theta d \phi ^2 ) \right)\,.
\]
Under assumption that cosmological constant is exactly zero the function  $a(t)$ can be found via parameter $\eta $ as follows, 
\[
a=a_*(1-\cos \eta ), \quad t= a_*(\eta -\sin \eta ), \quad a_*=\frac{2GM}{3\pi} ,
\]
which implies for the function $a=a(t)$ a periodic dependence on time. Although it is not known whether the universe actually has  exactly periodic behavior suggested by this 
solution of the Einstein equations, we find important to investigate {\it influence of the periodic behavior of metric tensor on the nonlinear waves and wave maps} 
propagating in this kind of spacetime. 
 \smallskip

Among  publications on the time periodic solutions of the Einstein's field equations we shall mention very few of them. The study of
the periodic solutions to the Einstein's field equations was initiated by Papapetrou in \cite{Papapetrou1957,Papapetrou-Treder,Papapetrou1962}.
Einstein's old question concerning the existence of solutions free of singularities, vanishing and becoming Euclidean at infinity is generalized in \cite{Papapetrou-Treder}
 for the case of periodic solutions free of singularities. 
 The result is negative; such solutions do not exist. 
  It is also shown by Gibbons and Stewart in \cite{Gibbons} that   any asymptotically flat spacetime which is empty and periodic in time is stationary in a neighborhood of null infinity.  
 
Dafermos~\cite{Dafermos} proved a theorem about the non-existence of spherically symmetric black-hole
spacetimes with time-periodicity outside the event horizon, other than Schwarzschild in the 
vacuum case and Reissner-Nordstr{\"o}m in the case of electromagnetic fields  and matter sources of a particular kind. 
This  result generalizes the so-called ``no-hair'' theorem from the static to the time-periodic case. That
 paper   addresses the issue of the existence of periodic solutions
in general relativity in a non-analytic setting, in particular, in a setting compatible
with the evolutionary hypothesis.
For the equations of evolution, time-periodic or stationary solutions often correspond to the
late time behavior of solutions for a large class of initial data. In the general theory of
relativity, time-periodic "black hole" solutions, if they exist, seem to provide reasonable
candidates for the final state of gravitational collapse. 
(See \cite{Dafermos,Kong_I,Kong_II}  for more references.)

 \smallskip

For the nonlinear scalar waves and the  scalar wave map type equations the influence of the periodic behavior of metric tensor on the nonlinear waves and wave maps 
 propagating in this kind of universe was found in \cite{YagJMAA2001,yagdjian_birk}.
It was discovered in those papers that, the generated by periodic coefficient parametric resonance interacting with the nonlinearity,
in general, leads to the blow up of the solution for arbitrary small initial data and for any dimension of the space of spatial variables. 

In order to establish similar result for the wave maps we first  prove (Lemma~\ref{L3.3}) that, in the case of $n=2$ 
the only stationary  wave maps  are the constant  wave maps. They evidently have vanishing integral (\ref{epsilon}). 
This is why we are interested in the 
small  perturbations of the constant  wave maps (with the periodic in time metric $g$).

It was proved in \cite[Lemma~2.1]{Krieger_arxiv2004} that for $l=2$, $n=2$ and Minkowski space with metric (\ref{MinMetr}), 
the image of the  spherically symmetric wave map with smooth initial
data  belongs to a bounded subset of ${\mathbb H}^2$. In particular, for the second component of the wave map we have
\begin{eqnarray} 
\label{lnu2} 
\| \ln u^2 \|_{L_t^\infty L_x^\infty}<\infty \,.   
\end{eqnarray} 
 
The main result of the present note is the following theorem that 
requires some technical assumption (Assumption ISIN, see Section~\ref{S5} for details) on the ordinary 
differential equation generated by the periodic metric. Without loss of generality, we can assume that
\begin{eqnarray}
\label{at0}
R'(0)=0\,.
\end{eqnarray}
\begin{theorem} 
Let $R=R(t)$ of (\ref{R-W_Metr}) be periodic, positive, non-constant, smooth function satisfying assumption ISIN, and let be (\ref{at0}). 
Denote  $({\mathbb H}_l^2,h)$
the Riemannian manifold,
that is, the  half-plane  ${\mathbb H}_l^2:=\{(u^1,u^2) \in {\mathbb R}^2 \,|\, u^2>0 \}$
equipped with the metric  tensor   $\dsp h_{AB}d u^A d u^B = \frac{1}{(u^2)^l}\left( (d u^1)^2+ (d u^2)^2\right)$. 
Then, the local solution obtained by small,  in the sense of the following integral 
\begin{eqnarray}
\label{epsilonG2}
\int_{\{0\}\times  {\mathbb R}^n}\sum_{   \gamma_0+ \gamma_1+\gamma_2+ \ldots  +\gamma_n =1,\ldots,G
\atop{\gamma _0 =0,1} }   \left( \left[ \frac{\partial^\gamma  u^1}{(u^2)^{l/2}} \right]^2+  \left[ \frac{\partial^\gamma u^2}{(u^2)^{l/2}} \right]^2\right) dx  ,
\end{eqnarray}
some smooth compactly supported spherically symmetric perturbations of the initial data of the constant wave map $(V,g) \longrightarrow  ({\mathbb H}_l^2,h)$,   
cannot be extended  globally in time if $l\in [0,2)$. 

If   $l=2$ and $n=2$, then   for some smooth compactly supported, with the arbitrarily small integral (\ref{epsilonG2}),  
spherically symmetric perturbations of the initial data of the constant wave map, 
there exist  positive numbers $c_0 $ and $\delta  $ such that for all $m \in {\mathbb N} $ either
\begin{eqnarray*} 
\inf_{x \in {\mathbb R}^2  } \ln u^2 (x ,m)   \leq  -c_0 e^{\delta m}
\quad or \quad \sup_{x \in {\mathbb R}^2  } \ln u^2 (x ,m)    \geq c_0 e^{\delta m}\,.  
\end{eqnarray*}
\end{theorem}

Thus, in the case of $l=2$ and $n=2$, due to the periodicity of the metric of the spacetime, the problem in the neighborhood 
of the stationary  spherically symmetric wave map is unstable 
in the following sense: for some smooth compactly supported, with the arbitrarily small integral (\ref{epsilonG2}),  
spherically symmetric perturbations of the initial data of the  stationary  spherically symmetric wave map, the estimate (\ref{lnu2}) 
does not hold.

\medskip

The proof of the theorem is based on the construction, for every given positive integer number $G$, of the wave map, 
which takes arbitrarily close to the constant initial data but the image $u(t_{bp},x_{bp}) $ 
of some point $x_{bp} \in {\mathbb R}^n$ in the finite time $t_{bp}>0$ appears  outside of the target  manifold  ${\mathbb H}_l^2$, that is, 
  $u^2(t_{bp},x_{bp})\le 0 $. In that sense the case of $l=0$, 
which implies a linear system for two components of $u(t ,x ) $, is not exeptional and the component $u^2 $
has a finite life span (becomes nonpositive in the finite time). 
\medskip

This note is organized as follows.  In Section~\ref{S2} we describe the families of geodesics of the target space.
This geodesics we use to construct wave maps via linear wave equation. 
In Section~\ref{S3}, we analyze the structure of the  wave map equations in the Minkowski and Robertson-Walker spacetimes. 
 Section~\ref{S4}    is devoted to the finding of  stationary wave maps.
 In Section~\ref{S5} in order to prove an exponential growth of   $L_x^\infty$-norm of the solutions of the scalar wave equation with periodic 
coefficients, we derive   explicit representations of some increasing solutions of 
the ordinary differential equations with periodic coefficients.
 In the final Section~\ref{S6} we complete
the proof of the main result.    
\medskip

\section{The Target Space. Geodesics}
\label{S2}

Let ${\mathbb H}_l^2:=\{(u^1,u^2)\,|\, u^2>0 \}$ be a half-plane equipped with the metric  
\[
\dsp h_{ij}\,du^i\, du^j = \frac{1}{(u^2)^l }\left( (du^1)^2+(du^2)^2\right),  
\]
where $l \in {\mathbb R}$ is a real number, while the metric tensor is 
\begin{eqnarray*}
\left( h_{AB} \right)_{A,B=1,2}=\frac{1}{(u^2)^l}\pmatrix{1 & 0\cr 0 & 1 } =\frac{1}{(u^2)^l}I \,.  
\end{eqnarray*}
If $l=2$ this is the   Poincar{\'e}  upper half-plane model. The upper half-plane  ${\mathbb H}$
endowed with the hyperbolic metric  with  $l=1$ is discussed in Ch.3~\cite{Struwe}.  The  Christoffel symbols are
\begin{eqnarray*}
\left( \Gamma _{A B }^1\right)_{A,B=1,2} 
= -  \frac{l}{2u^2}\pmatrix{0 & 1\cr  1 &0 },\quad 
\left( \Gamma _{A B  }^2 \right)_{A,B=1,2}
= -  \frac{l}{2u^2}\pmatrix{-1 & 0\cr 0 & 1 } .
\end{eqnarray*}
The Gaussian curvature is 
\[
K_l =   -\frac{l}{2}(u^2)^{l-2} .
\]
Hence, in this family of Riemannian manifolds only the Poincar{\'e}  upper half-plane model and the Euclidean half-plane have  constant curvature $-1 $ and $0$, respectively.

The geodesics $(u^1,u^2) =(u^1(s),u^2(s))$ satisfy equations 
\[
\frac{d^2 u^1}{ds^2} -  \frac{l}{u^2} \frac{d  u^1 }{ds }\frac{d  u^2  }{ds }=0,\quad \frac{d^2 u^2}{ds^2}+ \frac{l}{2u^2} \frac{d  u^1 }{ds }\frac{d  u^1  }{ds }
- \frac{l}{2u^2} \frac{d  u^2 }{ds }\frac{d  u^2  }{ds }=0\,,
\]
where $s$ is a natural parameter. If we denote $\dsp V_1:=\frac{d  u^1 }{ds }$ and $\dsp V_2:=u^2>0$, then
\[
\frac{d  V_1}{ds } -  \frac{l}{V_2} V_1  \frac{d  V_2  }{ds }=0,\quad \frac{d^2 V_2}{ds^2}+ \frac{l}{2V_2} (V_1)^2 
- \frac{l}{2V_2} \left( \frac{d  V_2 }{ds } \right)^2 =0\,.
\]
From the first equation we have 
\[
 V_1= C \left(  V_2   \right)^l\,,\qquad C=const \in {\mathbb R}.
\]
By plugging this into the unite speed relation
\[ 
1=\frac{1}{(V_2)^l}(V_1)^2+ \frac{1}{(V_2)^l}((V_2)_s)^2\,,
\]
we obtain
\begin{eqnarray*}
&  &
(V_2)^l- C^2\left(  V_2   \right)^{2l}= ((V_2)_s)^2
\end{eqnarray*}
with the constraint 
\[ 
C^2\left(  V_2 (t)   \right)^{ l} \leq 1  \,.
\]
 Thus, the function $V_2=V_2(s)>0 $ solves the following equation
\begin{eqnarray*}
&  &
\frac{dV_2}{ds}  =  \pm \sqrt{ (V_2)^l- C^2\left(  V_2   \right)^{2l}}\,.
\end{eqnarray*} 
Consequently, the  system 
\begin{eqnarray}
\label{geo_system}
&  &
\cases{
\dsp \frac{du^2}{ds}  =  \pm \sqrt{ (u^2)^l- C^2\left(  u^2   \right)^{2l}}\cr 
{}\cr
\dsp \frac{du^1}{ds} = C \left(  u^2   \right)^l} 
\end{eqnarray}
has the solution 
\[
u^1(t)= u^1(0) , \qquad u^2(t)=  \left[(u^2(0))^{\frac{2-l}{2}} \pm \frac{2-l}{2} s\right]^{\frac{2}{2-l}}  \quad \mbox{\rm if} \quad C=0,\quad l\not=2,
\]
or
\[
u^1(s)= u^1(0) , \qquad \left( \left[u^2(s)\right]^{\frac{2-l}{2}}-\left[u^2(0) \right]^{\frac{2-l}{2}}\right)^2=    
 \left(\frac{2-l}{2} s\right)^2  \quad \mbox{\rm if} \quad C=0,\quad l\not=2,
\] 
and
\[
u^1(s)= u^1(0) , \qquad u^2(s)=   u^2(0) e^{s}  \quad \mbox{\rm if} \quad C=0,\quad l =2,
\]
that is a vertical open half-line in the positive half-plane.   For $C\not=0$ one can eliminate the variable $s$ and rewrite  the system  as 
a single equation
\begin{eqnarray*}
&  & 
\dsp \frac{du^2}{du^1 }  =  \frac{\pm \sqrt{ (u^2)^l- C^2\left(  u^2   \right)^{2l}}}{C \left(  u^2   \right)^l}  \,.
\end{eqnarray*}
Then, one can integrate it
\begin{eqnarray*}
&  &
\pm \int \dsp\frac{C \left(  u^2   \right)^{l/2}}{\sqrt{ 1- C^2\left(  u^2   \right)^{l}}}  du^2 = \int  du^1\,.
\end{eqnarray*}
If we denote $u:= u^1 $  and $ v:=u^2$, then  
\begin{eqnarray*}
&  &
\pm \int_0^v \dsp\frac{C   x  ^{l/2}}{\sqrt{ 1- C^2 x ^{l}}}  dx = u- C_1 \,.
\end{eqnarray*}
We make change of dummy variable $x^l=v^lt$ with $lx^{l-1}dx=v^ldt $ in the  integral,
\begin{eqnarray*} 
 \int_0^v \dsp\frac{C   x  ^{l/2}}{\sqrt{ 1- C^2 x ^{l}}} \, dx  
& =  &
   v^{\frac{2+l}{2 } }  \frac{C}{l}\int_0^{1} \dsp\frac{ t^{\frac{2-l}{2l} }   }{\sqrt{ 1- (C^2v^l) t }} \, dt 
\,.  
\end{eqnarray*}
To evaluate the last integral we use formula (10)~Sec.2.1.3~\cite{B-E} with
\[
a=\frac{1}{2}, \quad b= \frac{2+l}{2l }, \quad c= 
\frac{ 3 }{2  }+   \frac{ 1}{ l }, \quad c-b-1=0, \quad z= C^2v^l, \quad 
\frac{\Gamma (b) \Gamma (c-b)}{\Gamma (c)} 
=  \frac{2l }{2+ l} \,,
\]
and obtain
\begin{eqnarray} 
\label{9}
u- C_1 = \pm \frac{2}{2+l } C v^{\frac{2+l}{2}}  F \left(\frac{2+l }{2 l },\frac{1}{2};  \frac{3}{2} +\frac{1}{l } ;C^2v^{l}\right) \,,
\end{eqnarray}
where $F\big(a, b;c; \zeta \big) $ is the hypergeometric function (See, e.g., \cite{B-E}.).
In fact,
\[
(u^1- C_1)^2 = \left( \frac{2}{2+l }\right)^2 C^2 (u^2)^{ {2+l} } \left[F \left(\frac{2+l }{2 l },\frac{1}{2};  \frac{3}{2} +\frac{1}{l } ;C^2(u^2)^{l}\right) \right]^2\,.
\]
In particular, for $l=1$ we obtain,
\begin{eqnarray*}
(u^1- C_1)^2 =  \frac{1}{C^4  } \left[-C \sqrt{  u^2  -C^2 (u^2)^2 } +\arcsin  \left(C \sqrt{ u^2 }\right) \right]^2\,,
\quad 0< u^2 \leq C^{-2}\,.
\end{eqnarray*}
For $l=2$ the geodesics are vertical lines and the  upper half-circles given by equation
\[
\left(u^1-d \right)^2+(u^2)^2=\frac{1}{C^2}\,, \qquad d= C_1+\frac{1}{C }\,.
\]
The   upper half-circles can be written   via the parameter $s$ as follows
\begin{eqnarray} 
\dsp u^1(s)= d+\frac{1}{C}\frac{e^{2s}-1}{e^{2s}+1}, \qquad
\dsp u^2(s)= \frac{2}{C} \frac{e^{s} }{e^{2s}+1}, \quad s \in {\mathbb R}\,.  
\end{eqnarray}
Thus, we have obtained the second family of geodesics.

 To find out the function $u^2= u^2(s)  =v(s) $ we use the first equation of the system (\ref{geo_system}) and, similar to the derivation of (\ref{9}), obtain
\begin{eqnarray*}
&  &
\frac{2 }{2-l } \sqrt{1-C^2 v^{  l}} \left\{ v^{1-\frac{l}{2}}+\frac{2}{2+l } C^2 v^{\frac{2+l}{2}}  F\left(1,1+\frac{1}{l};\frac{3}{2}+\frac{1}{l};C^2 v^l\right)\right\}\\
& - &
\frac{2 }{2-l } \sqrt{1-C^2 v^{  l}(0)} \left\{ v^{1-\frac{l}{2}}(0)+\frac{2}{2+l } C^2 v^{\frac{2+l}{2}}(0)  F\left(1,1+\frac{1}{l};\frac{3}{2}+\frac{1}{l};C^2 v^l(0)\right)\right\}
=s\,,
\end{eqnarray*}
which for the positive $u^2(0) $ with $C^2 (u^2(0))^{  l}  \leq 1$ defines an implicit function $u^2=u^2(s) $. Then one can use (\ref{9}) to find out the function $u^1=u^1(s) $.

The composition of a geodesic with any real-valued
solution $\psi=\psi (x,t) $ of the free wave equation 
\[
\frac{1}{\sqrt{|g|}}\frac{\partial }{\partial x^i}\left( \sqrt{|g|} g^{ik} \frac{\partial \psi }{\partial x^k} \right) = 0
\] 
generates the wave map $u = u_\gamma $ (see \cite{Shatah}). Such special solutions  have been used in \cite{Piero-Georgiev} to prove, for instance, 
that solution map is not locally Lipschitz continuous.  
Let $ \gamma (s)$ be an arbitrary curve with values in ${\mathbb H}^2_l$, which is written in the local coordinates by $ \gamma (s)= (\gamma_1 (s),\gamma_2 (s) )$ and 
 let $v(x,t)$ be an arbitrary real-valued function, then 
the function $u(x,t)= \gamma (v(x,t))$ solves the wave map equation, as soon as $v(x,t) $ solves the wave equation and $\gamma (s) $ is a geodesic curve. 
Thus, we have proved the following statement.
\begin{theorem}
\label{L1.1}
\label{T1.1}
Let function $\varphi =\varphi (x,t) $ be a solution of the covariant  wave equation in $(V,g)$
\begin{eqnarray*} 
&  & 
g ^{\alpha \beta } \left(  \partial_{\alpha\beta}^2  \varphi  - \Gamma _{\alpha  \beta }^\lambda  \partial_\lambda   \varphi 
   \right)=0\,.
\end{eqnarray*}
Then the following pairs $(u^1,u^2) $ of functions,
\begin{eqnarray*}
&  &
u^1(x,t)=  C_1,\qquad u^2 (x,t)=  \left[(C_2)^{\frac{2-l}{2}} \pm \frac{2-l}{2} \varphi (x,t)\right]^{\frac{2}{2-l}}, \quad t \in {\mathbb R}, \,\, x \in {\mathbb R}^n,  \quad 
 l \not=2,\\
&  &
\dsp u^1(x,t)= C_1, \qquad
\dsp u^2(x,t)=   e^{\varphi (x,t)}  , \quad t \in {\mathbb R}, \,\, x \in {\mathbb R}^n,\quad   l =2,\,\\
&  &
\dsp u^1(x,t)= d+\frac{1}{C}\frac{e^{2\varphi (x,t)}-1}{e^{2\varphi (x,t)}+1}, \qquad
\dsp u^2(x,t)= \frac{2}{C} \frac{e^{\varphi (x,t)} }{e^{2\varphi (x,t)}+1}, \quad t \in {\mathbb R},\,\, x \in {\mathbb R}^n,\quad   l =2\,,   
\end{eqnarray*}
and the function $(u^1(x,t), u^2(x,t))=(u^1(x,t), v(x,t)) $ defined by 
\begin{eqnarray*}
&  &
u^1(x,t)=   C_1  \pm \frac{2}{2+l } C v (x,t)^{\frac{2+l}{2}}  F \left(\frac{2+l }{2 l },\frac{1}{2}; \left(\frac{3}{2} +\frac{1}{l }\right);C^2v (x,t) ^{l}\right), \\
&  &
\frac{2 }{2-l } \sqrt{1-C^2 v^{  l}} \left\{ v^{1-\frac{l}{2}}+\frac{2}{2+l } C^2 v^{\frac{2+l}{2}}  F\left(1,1+\frac{1}{l};\frac{3}{2}+\frac{1}{l};C^2 v^l\right)\right\}\\
& - &
\frac{2 }{2-l } \sqrt{1-C^2 v^{  l}(0)} \left\{ v^{1-\frac{l}{2}}(0)+\frac{2}{2+l } C^2 v^{\frac{2+l}{2}}(0)  F\left(1,1+\frac{1}{l};\frac{3}{2}+\frac{1}{l};C^2 v^l(0)\right)\right\}
=\varphi (x,t)
\end{eqnarray*}
are the wave maps $u: V\longrightarrow  {\mathbb H}^2_l$.
\end{theorem}

\section{The Wave Maps in the Robertson-Walker Spacetimes}
\label{S3}

In the metric 
\begin{eqnarray*} 
&  &
g   = -dt^2+R^2(t) \sigma, 
\end{eqnarray*}
in the local coordinates the wave map equation (\ref{WME}) 
reads
\[
-  \partial_t^2  u^A - nR^{-1} R' \partial_t   u^A+ R^{-2} \Delta _\gamma   u^A  
 - \Gamma _{BC }^A  \partial_t    u^B  \partial_t     u^C 
+ R^{-2}\sigma  ^{ij} \Gamma _{BC }^A  \partial_i    u^B  \partial_j     u^C =0\,,
\]
where we denote by  
\[
 \Delta _\gamma = \sigma  ^{ij }   \partial_{ij}^2  
 -   \sigma  ^{ij }\gamma _{ij }^k  \partial_k 
\]
  the Laplace operator on the manifold with metric $\sigma  ^{ij} $. The wave map solves the system
\[ 
\cases{\dsp 
-  \partial_t^2  u^1 - nR^{-1} R' \partial_t   u^1+ R^{-2} \Delta _\gamma   u^1  
 + \frac{l}{ u^2 }  \partial_t    u^1  \partial_t     u^2 
- R^{-2}\sigma  ^{ij} \frac{l}{ u^2 } \partial_i    u^1  \partial_j     u^2 =0,\cr
{}\cr 
\dsp -  \partial_t^2  u^2 - nR^{-1} R' \partial_t   u^2+ R^{-2} \Delta _\gamma   u^2  
 -  \frac{l}{2u^2 }( \partial_t    u^1  \partial_t     u^1 - \partial_t    u^2  \partial_t     u^2)\cr
\dsp \hspace{6cm} + R^{-2}\sigma  ^{ij}  \frac{l}{2u^2 } (  \partial_i    u^1  \partial_j     u^1 -  \partial_i    u^2  \partial_j     u^2)=0.}
\]
If $ \sigma  ^{ij}= \delta   ^{ij}$ then we obtain 
\begin{equation}
\label{system7}
\cases{\dsp
  \partial_t^2  u^1 + nR^{-1} R' \partial_t   u^1 - R^{-2} \Delta     u^1  
 - \frac{l}{ u^2 } \partial_t    u^1  \partial_t     u^2 
+ R^{-2} \frac{l}{ u^2 } \nabla_x    u^1 \cdot \nabla_x    u^2 =0,\cr 
{}\cr
\dsp \partial_t^2  u^2 + nR^{-1} R' \partial_t   u^2- R^{-2} \Delta     u^2  
 +   \frac{l}{2u^2 }( \partial_t    u^1  \partial_t     u^1 - \partial_t    u^2  \partial_t     u^2)\cr
\dsp \hspace{6cm} 
- R^{-2}  \frac{l}{2u^2 } (  |\nabla_x    u^1  |^2 -  |\nabla_x    u^2  |^2)=0.}
\end{equation}
From now on we say that the wave map $(u^1,u^2)$ has a finite life span (blows up in the finite time) if for some point $(x_0,t_0)$ with $t_0>0$ 
the image 
$(u^1(x_0,t_0),u^2(x_0,t_0))$  appears  outside of the target  manifold ${\mathbb H}^2_l $, that is either $ u^2(x_0,t_0) = 0$ or $ u^2(x_0,t_0) = \infty$.

To reveal the blowup mechanism we use the first wave map of Theorem~\ref{T1.1}, that is, we set in the last system  $u^1\equiv C_1 =const$. 
This corresponds to the choice of the geodesic
$\gamma (s) $, which is the open vertical half-line.  The second equation of the system reads
\begin{eqnarray*}
&  &
\dsp \partial_t^2  u^2 + nR^{-1} R' \partial_t   u^2- R^{-2} \Delta     u^2  
 -   \frac{l}{2u^2 } \left(  (\partial_t    u^2)^2 
- R^{-2}      |\nabla_x    u^2  |^2  \right)=0\,.
\end{eqnarray*}
Consider the case of $l\not=2$. Then we can rewrite the last equation as follows:
\begin{eqnarray*}
\dsp \partial_t^2  u^2 + nR^{-1} R' \partial_t   u^2- R^{-2} \Delta     u^2  
 -  \frac{\mu -1}{\mu u^2 }\left(  (\partial_t    u^2 )^2 
- R^{-2}   (   |\nabla_x    u^2  |^2 \right)=0\,.
\end{eqnarray*}
Here 
\[
 \frac{l}{2 }= \frac{\mu -1}{\mu   } \,, \qquad \mu = \frac{2}{2 -l} \not=0\,.
\]
In particular, if  $l=1$, then $\mu =2$, while for $l=4$ we obtain $\mu =-1$.
If we denote $u:=u^1 $ and introduce a new unknown function $v$ by means of equation
\[
u^2= \left(\frac{\alpha }{\mu }v+\beta \right)^\mu >0,
\]
with parameter $\alpha \in {\mathbb R} $ and the constants $\beta  $, $l/2=(\mu -1)/\mu $, then from the equations (\ref{system7}) we obtain
\begin{eqnarray*}
\cases{\dsp
   u_{tt} + nR^{-1} R'   u_t - R^{-2} \Delta     u  
 -  l \alpha   \left(\frac{\alpha  }{\mu }v+\beta   \right)^{ -1} \left[      u_t   v_t  
- R^{-2}   \nabla_x    u  \cdot    \nabla_x v \right] =0,\cr 
{}\cr
\dsp       v_{tt}     + nR^{-1} R'    v_t   - R^{-2}     \Delta    v     +   \frac{l}{2\alpha   }\left(\frac{\alpha  }{\mu }v+\beta  \right)^{1-2\mu}\left[ (     u_t )^2  
- R^{-2}     |\nabla_x    u  |^2  \right]=0.}
\end{eqnarray*}
For the case of $l\not= 2$ and $u^1\equiv const $ we obtain that the function 
$v$ solves the following linear equation
\begin{eqnarray}
\label{LinEQ}
&  &
\partial_t^2  v + nR^{-1} R' \partial_t   v- R^{-2} \Delta    v  =0\,.
\end{eqnarray}
If we prove that at some point $(x_0,t_0)$ with $t_0>0$  the function $v(x,t) $ takes value $-\mu \beta  /\alpha $, that is, $v(x_0,t_0)= -\mu \beta  /\alpha  $, 
then at that point the wave map  $(u^1,u^2)$ blows up. 
Indeed, it is evident from the definition of $v$ that if $ \mu <0$ then  at that point $u^2(x_0,t_0)=\infty$. 
If   $\mu >0$, then at that point  we obtain $u^2(x_0,t_0)=0 $, and the image of the point $ (x_0,t_0)  $   is outside of the manifold
${\mathbb H}^2_l $. The last case is regarded as blow up as well. Note that the half-lines $u^1=const$, $u^2>0$, 
are geodesic in the target manifold ${\mathbb H}^2_l $.

\bigskip

\noindent
 In the case of the  Poincar{\'e}  upper half-plane model ${\mathbb H}^2_l $ we have  $l=2$ and 
\begin{eqnarray*}
\left( h_{AB} \right)_{A,B=1,2}=\frac{1}{(u^2)^2}\pmatrix{1 & 0\cr 0 & 1 } =\frac{1}{(u^2)^2}I,\qquad 
\left( h^{AB} \right)_{A,B=1,2}= (u^2)^2 \pmatrix{1 & 0\cr 0 & 1 } = (u^2)^2 I,
\end{eqnarray*}
\begin{eqnarray*}
\left( \Gamma _{A B }^1 \right)_{A,B=1,2}= -  \frac{1}{u^2}\pmatrix{0 & 1\cr  1 &0 } ,\qquad 
\left( \Gamma _{A B }^2 \right)_{A,B=1,2} =  -  \frac{1}{u^2}\pmatrix{-1 & 0\cr 0 & 1 }, 
\end{eqnarray*}
and the corresponding wave map equation is  
\[ 
\cases{ \dsp \partial_t^2  u^1 + nR^{-1} R' \partial_t   u^1- R^{-2} \Delta    u^1  
 - \frac{2}{ u^2}  (\partial_t    u^1 )( \partial_t     u^2 )
+ R^{-2}  \frac{2}{ u^2} \nabla_x   u^1 \cdot  \nabla_x     u^2 =0,\cr
{}\cr 
\dsp \partial_t^2  u^2 + nR^{-1} R' \partial_t   u^2- R^{-2} \Delta    u^2  
 +   \frac{1}{u^2}( (\partial_t    u^1 )^2 - (\partial_t    u^2 )^2) - R^{-2}  \frac{1}{u^2} (  |\nabla_x    u^1  |^2 -  |\nabla_x    u^2  |^2)=0.}
\]
If we set
\[ 
u:=u^1,\quad u^2= e^v >0\,,
\]
then we arrive at the system
\begin{eqnarray} 
\label{starsystem}
\cases{ \dsp    u_{tt} + nR^{-1} R'     u_t - R^{-2} \Delta _\gamma   u  
 - 2       u_t    v_t  
+ 2R^{-2}    \nabla_x   u \cdot      \nabla_x    v=0,\cr
{}\cr 
\dsp     v_{tt}   + nR^{-1} R' v_t   - R^{-2}    \Delta _\gamma  v    +    e^{-2v} \left[ (   u_t )^2
- R^{-2}    |\nabla_x    u  |^2  \right] =0.}
\end{eqnarray}
If we consider case with $u\equiv const$, then the second equation  implies  that 
$v$ solves linear equation (\ref{LinEQ}) 
which has a global solution. Thus, if the initial data for $u^1  $ are the constants, then a global solution of the non-linear equation (wave map) exists. On the other hand, if 
$u^2=const$ then the first equation is linear, and, consequently, it has a global solution.  These two families of lines form  semi-geodesic parametrization of the target manifold.

\medskip

\section{Stationary Solutions}
\label{S4}

We look for the stationary solutions of the wave map system
\begin{eqnarray}
\label{system_wm}
\cases{\dsp
  \partial_t^2  u^1 + nR^{-1} R' \partial_t   u^1 - R^{-2} \Delta _\gamma   u^1  
 - \frac{l}{ u^2 } \partial_t    u^1  \partial_t     u^2 
+ R^{-2} \frac{l}{ u^2 } \nabla_x    u^1 \cdot \nabla_x    u^2 =0,\cr 
{}\cr
\dsp \partial_t^2  u^2 + nR^{-1} R' \partial_t   u^2- R^{-2} \Delta _\gamma   u^2  
 +   \frac{l}{2u^2 }( \partial_t    u^1  \partial_t     u^1 - \partial_t    u^2  \partial_t     u^2) - R^{-2}  \frac{l}{2u^2 } (  |\nabla_x    u^1  |^2 -  |\nabla_x    u^2  |^2)=0 }
\end{eqnarray} 
with the positive function $u^2 (t,x)>0 $. The stationary wave map solves the following system of quasilinear elliptic equations
\begin{eqnarray}
\label{elliptic}
\cases{\dsp
\Delta _\gamma   u^1  
-  \frac{l}{ u^2 } \nabla_x    u^1 \cdot \nabla_x    u^2 =0,\cr 
{}\cr
\dsp   \Delta _\gamma   u^2  
+ \frac{l}{2u^2 } (  |\nabla_x    u^1  |^2 -  |\nabla_x    u^2  |^2)=0.}
\end{eqnarray} 
The second unknown function is assumed to be positive
that allows to invoke the Liouville theorem for the superharmonic functions. 
The following statement is evident. 

\begin{lemma}
\label{Levident}
For every $ l \in [0,2]$ the set of stationary solutions for wave map equation is independent of the choice of the function $R=R(t)$. 
\end{lemma}
Thus for the case of $l=2$ and $n=2$ we can appeal to the next lemma, which is due to \cite{Krieger_arxiv2004}.
\begin{lemma} \mbox{\rm (\cite{Krieger_arxiv2004})} 
\label{LKrieger_2006}
The image of the wave map belongs to a bounded subset of ${\mathbb H}^2$. More precisely, we have
\begin{eqnarray} 
\label{lemma_kreiger}
\| \ln u^2 \|_{L_t^\infty L_x^\infty}<\infty ,
\qquad \| \frac{u^1}{ u^2} \|_{L_t^\infty L_x^\infty}<\infty\,.  
\end{eqnarray}
The bounds depend (at most) on the size of the support as well as some norm $\|   u[0] \|_{H^{1+\delta }} $, $\delta >0 $.
\end{lemma}

\begin{lemma}
\label{L3.3}
Suppose that $n=2$ and $\dsp 0\leq l< 2  $. The only stationary solutions of the system  (\ref{system_wm}) with the finite integral
\begin{eqnarray} 
\label{norm_kreiger}
\int_{{\mathbb R}^2}\sum_{\alpha  =1}^2  \left( \left[ \frac{\partial_\alpha   u^1}{(u^2)^{l/2}} \right]^2+  \left[ \frac{\partial_\alpha  u^2}{(u^2)^{l/2}} \right]^2\right) dx 
\end{eqnarray}
and with the positive $u^2 (t,x)>0 $ are the constant solutions, $ u^1 (t,x) \equiv c_1 $,  $u^2 (t,x)\equiv c_2 >0 $.
If $l=2$, then  the only stationary spherically  symmetric wave maps with target ${\mathbb H}^2$
and finite integral (\ref{norm_kreiger}) are the constant wave maps.
\end{lemma}
\medskip

\noindent
{\bf Proof.} 
The stationary solution of (\ref{system_wm}), that is solution independent of $t$, $ u^1(t,x)= u^1(x)$,  $u^2(t,x)= u^2(x) $, solves the system (\ref{elliptic}). 
In the case of $l=0$ we have two harmonic in  ${\mathbb R}^2 $ functions with the finite integral (\ref{norm_kreiger}). The Liouville theorem (see, e.g. \cite[Theorem II]{Serrin-Zou})
for the positive function $u^2(x) $ 
implies  $u^2(x) =const>0$. Then, the harmonic functions $\partial_\alpha u^1 (x) $, $\alpha =1,2 $, belong to $ L^2({\mathbb R}^2)$ only if  $u^1(x) =const$. 

For the case of $ l \in (0,2) $ we set
\begin{eqnarray*}
u:=u^1,\qquad u^2=v^\mu , 
\end{eqnarray*}
where   $\mu =2/(2-l) $, and, consequently,
the system (\ref{elliptic}) reads 
\[
\cases{\dsp
\Delta    u  
-  \frac{l\mu}{ v } \nabla_x    u \cdot    \nabla_x v =0,\cr 
{}\cr
\dsp   \mu v^{\mu -1} \Delta   v+ \mu (\mu -1) v^{\mu -2} |\nabla_x v|^2  
+ \frac{l}{2v^\mu } \left(  |\nabla_x    u |^2 -  \mu^2v^{2\mu -2}|  \nabla_x v  |^2\right)=0.}
\]
Then we obtain 
\[
\cases{\dsp
\Delta   u  
-  \frac{ 2(\mu-1)}{ v } \nabla_x    u \cdot    \nabla_x v =0,\cr 
{}\cr
\dsp    \Delta   v  
+ \frac{\mu -1}{\mu^2   }  |\nabla_x    u |^2 v^{1 -2\mu}=0.}
\]
For the case of   $\dsp 0 < l<2  $ we have  $\dsp  \mu -1  > 0  $, and, consequently, the second equation of the last system implies
\[
\dsp    \Delta   v  
=- \frac{\mu -1}{\mu^2   }  |\nabla_x    u |^2 v^{1 -2\mu} \leq  0\,.
\]
Thus, the function $v=v(x) $ is superharmonic. According to the properties of superharmonic functions  (see, e.g., \cite{Serrin-Zou}), for the positive solution $v=v (x) $ the last equation implies  $v\equiv const>0 $. 
Then, the first equation of the system implies $ u= u^1 = u^1 (x)$ is harmonic in ${\mathbb R}^n $ function. 
The harmonic function has a finite integral (\ref{norm_kreiger}), that is, the derivatives of this function belong  to $L^2( {\mathbb R}^n ) $, only if  $ u^1 (x) \equiv const $. 

For $l=2$ we set $u:=u^1$, $u^2=e^v $ and then use (\ref{starsystem}). Hence we obtain   the stationary solution of the system
\begin{eqnarray*}
\cases{ \dsp   \Delta    u  - 
2\nabla_x   u \cdot      \nabla_x    v=0,\cr
{}\cr 
\dsp       \Delta   v    +   \frac{1}{ e^{2v}}|\nabla_x    u  |^2  =0\,.}
\end{eqnarray*}
Due to Lemma~\ref{Levident} and Lemma~\ref{LKrieger_2006}, we have $v(x) \geq \ln c $, where $c>0$. The non-negative function $w(x):=v(x) - \ln c \geq 0 $ solves equation
\begin{eqnarray*}
\dsp       \Delta _\gamma  w   +   \frac{1}{ ce^{2w}}|\nabla_x    u^1  |^2  =0\,. 
\end{eqnarray*}
Then we repeat the arguments have been used in the proof of the case of $0< l<2 $.
Lemma is proved. \hfill $\square$
\medskip

For $n >2$ there is a non-constant  bounded superharmonic in ${\mathbb R}^n $ function. (See, e.g., \cite{Gilbarg-Trudinger}.) 
Moreover,  the lemma implies the following statement. 
\begin{corollary}
Suppose that $\dsp 0\leq l< 2  $. 
If the stationary wave map $(V,g) \longrightarrow  ({\mathbb H}_l^2,h)$  with the finite integral (\ref{norm_kreiger})
is non-constant, then $n >2$.

\end{corollary} 
\medskip

This lemma motivates our interest to the perturbations, which have the finite integrals of type (\ref{norm_kreiger}), of the constant solutions.

\section{Parametric Resonance in ODE}
\label{S5}

Next we make the partial Liouville  transformation which eliminates the first derivative $v_t $ in (\ref{LinEQ}) with $\beta (x)=const $. 
 More precisely, we set   
\begin{eqnarray*}
&  &
v= R^{-\frac{n}{2}}w\,,
\end{eqnarray*}
then
\[
 v_{tt}- R^{-2}(t) \Delta v+ n R'R^{-1} v_t
= R^{-\frac{n}{2}}\Bigg[ w_{tt}- R^{-2}(t) \Delta w  +\left\{ -\frac{n}{2}R'{}'R^{-1}+ \frac{n}{2}\left(1-\frac{n}{2}  \right) (R')^2R^{-2} \right\} w  \Bigg]\,.
\]
Thus, we have to study the following linear hyperbolic equation 
\begin{eqnarray*}
 w_{tt}- R^{-2}(t) \Delta w  +\left\{ -\frac{n}{2}R'{}'R^{-1}+ \frac{n}{2}\left(1-\frac{n}{2}  \right) (R')^2R^{-2} \right\} w 
=0 
\end{eqnarray*}
with  the 1-periodic  positive smooth function $R=R(t)$.

To this end we are going to apply the Floquet-Lyapunov theory for the ordinary differential equation with the periodic coefficients. 
Consider the ordinary differential equation:
\begin{eqnarray*}
 W_{tt}+\left\{ \lambda   R^{-2}(t) -\frac{n}{2}R'{}'(t)R^{-1}(t)-\frac{n}{2}\left(\frac{n}{2}- 1 \right) (R'(t))^2R^{-2} (t)\right\} W 
=0 
\end{eqnarray*}
with the periodic positive smooth non-constant function $R=R(t)$ and parameter $\lambda \in {\mathbb R} $.

It is more convenient  to rewrite this equation by means of the new positive periodic function   
\[
\alpha (t)=  R^{-2}(t), \qquad  R (t)=  (\alpha(t))^{-1/2},
\]
then 
\begin{eqnarray*}
 W_{tt}+\left\{ \lambda   \alpha (t) - \frac{n}{4}\left[ \frac{3}{2} \left( \frac{\alpha '(t)}{\alpha (t)}\right)^2- \frac{\alpha'{}'  (t)}{\alpha (t)} \right]  
- \frac{n}{8}\left(\frac{n}{2}- 1 \right) \left( \frac{\alpha '(t)}{\alpha (t)}\right)^2  \right\} W 
=0 \,.
\end{eqnarray*}
Consider now the equation
\begin{eqnarray}
\label{Eastham1.2.1}
y_{tt}(t) + \left( \lambda    \alpha (t)-q(t) \right) y(t)=0
\end{eqnarray}
with  the periodic coefficients $\alpha (t)
  =   
R^{-2}(t) $ and  
\begin{eqnarray*} 
 q(t)
& =  & 
\frac{n}{4}\left[ \frac{3}{2} \left( \frac{\alpha '(t)}{\alpha (t)}\right)^2- \frac{\alpha'{}'  (t)}{\alpha (t)} \right]  
+ \frac{n}{8}\left(\frac{n}{2}- 1 \right) \left( \frac{\alpha '(t)}{\alpha (t)}\right)^2 \,.
\end{eqnarray*}
The first part of the last expression is the so-called Schwarz derivative for the antiderivative of $\alpha (t)$. For equation (\ref{Eastham1.2.1}) the spectrum of the  eigenvalue problem
\[
y(0)=y(1)=0
\]
is discrete.

The equation (\ref{Eastham1.2.1})  can be written also as a system of differential equations for the
vector-valued function $x(t)={}^t(w_t,w)$:
\[
\displaystyle{\frac{d}{dt}} x(t) = A(t)
 x(t)\,,\qquad  \mbox{\rm where}
\quad A(t) := \left(\begin{array}{lll}
0  & - \lambda    \alpha (t)+q(t) \\
1  & 0 
\end{array} \right)\,.
\]

Let the  matrix-valued function\, $X_\lambda (t,t_0)$,\, depending on \,$\lambda $,\,
be a solution of the 
Cauchy problem
\begin{equation}
\label{system}
\displaystyle{\frac{d}{dt}} X = A(t)
 X\,,\qquad 
X(t_0,t_0) = 
\left(\begin{array}{lll}
1  & 0 \\
0  & 1 
\end{array} \right)\,.
\end{equation}
Thus, \, $X_\lambda (t,t_0)$\, gives a fundamental solution  to the 
equation (\ref{Eastham1.2.1}).  
In what follows we often omit subindex $\lambda $ of $X_\lambda  (t,t_0)$. 
The Liouville formula 
\begin{eqnarray*}
W(t) = W(t_0) \exp \left( \int_{t_0}^t S(\tau) d\tau \right),\\ 
\end{eqnarray*}
 where $
W(t):= \det X(t,t_0)$, \,
$S(t):= \sum_{k=1}^2 A_{kk}(t)$ 
with $S(t) \equiv 0$ guarantees the existence of the inverse matrix
$X_\lambda  (t,t_0)^{-1}$. 
 For the  matrix $ X(1,0) $ we
will use a notation
\[
X_\lambda (1,0) =  \left(\begin{array}{lll}
 b_{11}  & b_{12}   \\
b_{21}      & b_{22}  
\end{array} \right)\,.
\]
This matrix is called a {\it monodromy matrix} 
and its eigenvalues are called {\it multipliers} of system (\ref{system}). Thus, the monodromy matrix is the value at 
$t=1$ (the ``end''of the period) of the fundamental matrix $X(t,0)$ 
defined by the initial condition $X(0,0)=I$ (i.e. the  matrizant), 
and the multipliers are the roots of the equation
\begin{eqnarray}
\label{ChE}
\det \left[ X(1,0) - \mu I \right] =0\,. 
\end{eqnarray} 
Due to Theorem~2.3.1~\cite{Eastham} there exist the   open instability intervals. We make the following 
\smallskip

\noindent  
{\bf Assumption ISIN:} {\it There exists the nonempty   open instability interval $\Lambda \subset (0,\infty)$ for equation (\ref{Eastham1.2.1}).} 
\smallskip

\noindent  
One can find in \cite{Eastham,Magnus_Winkler} detailed description of functions $\alpha =\alpha (t) $ and $q=q(t) $ satisfying this condition.
For instance, in Theorem~4.4.1~\cite{Eastham} one can find asymptotic formula, which allows to estimate the length of the instability intervals 
of the equation obtained from (\ref{Eastham1.2.1}) by Liouville transformation. 
Then, according to next lemma 
 one can find in 
the instability interval $\Lambda $  a 
number\,  $\lambda $\,  such  that a non-diagonal element of the monodromy matrix 
does not vanish. Moreover, this property is stable under small
perturbations of $\lambda $.

\begin{lemma} \mbox{\rm  (\cite{KY_Torino}-\cite{yagdjian_birk})}
\label{Lb21}
Let  \,  $R(t)$\, be defined on \, ${\mathbb R}$\, non-constant,  positive, 
smooth function which is  $1$-periodic.  Then there exists an open 
subset $\Lambda  ^0 \subset \Lambda  $
such that  $b_{21} \neq 0$ for all $\lambda  \in \Lambda ^0$.
\end{lemma}  
\medskip

Next we use the
 periodicity of \,$b = b(t)$\, and the eigenvalues \, $\mu_0>1$,\,  
$\mu_0^{-1}<1$
of the matrix 
$\,X_\lambda (1,0) \,$   
to construct solutions of (\ref{Eastham1.2.1}) with prescribed values on a discrete
set of time. The eigenvalues of matrix $X_\lambda (1,0)$ are $\mu_0$ and 
$\mu_0^{-1}$ with $ b_{11} + b_{22} = \mu_0 + \mu_0^{-1}$. 
Hence  $(b_{11} - \mu_0) + (b_{22} -\mu_0) = - \mu_0 + \mu_0^{-1}$
implies $|b_{11} - \mu_0| + |b_{22} - \mu_0| \geq 
|(b_{11} - \mu_0) + (b_{22} -\mu_0)| = |\mu_0 - \mu_0^{-1}| > 0$.
This leads to
\[
\max \{ |b_{11} - \mu_0| , |b_{22} - \mu_0|\} 
\geq \frac{1}{2} |\mu_0 - \mu_0^{-1}| > 0\,.
\]
Without loss of generality we can suppose
\[
|b_{11} - \mu_0| 
\geq \frac{1}{2} |\mu_0 - \mu_0^{-1}| > 0\,, \quad |b_{22} - \mu_0^{-1}| 
\geq \frac{1}{2} |\mu_0 - \mu_0^{-1}| > 0\,,
\]
because of $b_{11} - \mu_0 = -(b_{22} - \mu_0^{-1})$. Further, 
\[
1- \frac{b_{21}}{\mu_0^{-1}-b_{22}}\frac{b_{12}}{\mu_0-b_{11}}
= (\mu_0 - \mu_0^{-1})\frac{1}{b_{22} - \mu_0^{-1}}\not= 0 \,.
\]

\begin{lemma}\mbox{\rm (\cite{yagdjian_birk})}
Let $W = W(t)$, $V = V(t)$ be   solutions of the equation   
\[
 w_{tt} + \left( \lambda    \alpha (t)-q(t) \right) w  = 0 
\]
 with the parameter $\lambda $ 
such that $b_{21} \not= 0$ and $b_{22}\not=\mu_0^{-1} $. Suppose then that $W = W(t)$   takes  the initial data
\[
 W(0) = 0\,,\quad
 W_t(0) =  1  \,,
\]
and that $V = V(t)$ takes the initial data 
\[
V(0) = 1\,,\quad
 V_t(0) =  0  \,.
\]
Then for every positive integer number \,$M \in {\mathbb N}$\,  one has
\begin{eqnarray*}
 W( M) 
& = &
 \frac{b_{21}}
{\mu_0-\mu_0^{-1}}( \mu_0^M - \mu_0^{-M} ) \,,\\ 
 V( M) 
& = &
 - \mu_0^M\frac{(b_{22}-\mu_0^{-1})}{(\mu_0-\mu_0^{-1})}
 +  \mu_0^{-M} \frac{b_{21}b_{12}}{(\mu_0-b_{11})(\mu_0-\mu_0^{-1})} .
\end{eqnarray*}
\end{lemma}

In order to prove parametric resonance phenomena in $L^q$ spaces with $q\geq 2$, we need 
the following 

\begin{theorem}
\label{TLinRes} 
Let  \,  $R(t)$\, be defined on \, ${\mathbb R}$\, non-constant,  positive, 
smooth function which is  $1$-periodic.  
Then  for every function $ \varphi \in C_0^\infty ({\mathbb R}^n)$, which is different from the identical zero,
 there are positive numbers  $C_\varphi $, $\delta_\varphi  $   such that for $q \in [2,\infty]$ the solutions $w(x,t) $ and $v(x,t) $ 
of the equation 
\begin{eqnarray}
\label{3.16}
 w_{tt}- R^{-2}(t) \Delta w  +\left\{ -\frac{n}{2}R'{}'R^{-1}+ \frac{n}{2}\left(1-\frac{n}{2}  \right) (R')^2R^{-2} \right\} w 
=0 
\end{eqnarray}
with the initial conditions 
\begin{eqnarray*}
&  &
w(x,0)= \varphi (x )  ,\qquad w_t(x,0)=0,\\
\label{24}
&  &
 v(x,0)= 0,\qquad v_t(x,0)=\varphi (x  )   
\end{eqnarray*}
satisfy
\begin{eqnarray*}
&  &
\|  w(x,m) \|_{L^q({\mathbb R}^n)}  \geq  
  C_\varphi   e^{ \delta_\varphi  m}  \qquad \forall m \in{\mathbb N}, \\
&  &
\| v(x,m) \|_{L^q({\mathbb R}^n)} \geq
 C_\varphi   e^{ \delta_\varphi  m}   \qquad \forall m \in{\mathbb N}. 
\end{eqnarray*}
The constant  $\delta_\varphi  $ depends on the diameter of supp\,$\varphi  $ only. 
\end{theorem} 
In fact, $\mu_0=\mu_0(\lambda ) $ , $ b_{ij}=b_{ij}(\lambda ) $, $ i,j=1,2 $,   depend on $\lambda $ as well.
\medskip

In fact,  the functions   $W$ and $V$ depend on $\lambda $ as well,   $W = W(t,\lambda )$, $V = V(t,\lambda )$.
\medskip

\noindent
{\bf Proof.} Let $[a,b] \subset \Lambda ^0$ and $ \varphi \in C_0^\infty ({\mathbb R})$. 
Then by Paley-Wiener theorem  
\begin{eqnarray*}
&  &
\widehat \varphi (\xi ) \not\equiv 0\qquad \mbox{on }    I:=\left[\frac{a}{n},\frac{b}{n}\right]^n, \quad \mbox{\rm where}\quad |\xi|^2 \in [a,b] \quad \mbox{\rm for all}\quad 
\xi \in I,
\end{eqnarray*}
and for the Fourier transforms of the solutions we have
\begin{eqnarray*}
&  &
w(x,t):= (2\pi )^{n/2} \int e^{ix\xi}W(t,\lambda )\widehat \varphi (\xi )\, d\xi, \\
&  &
v(x,t):= (2\pi )^{n/2} \int e^{ix\xi}V(t,\lambda )\widehat \varphi (\xi )\, d\xi\,.
\end{eqnarray*}
Consequently, with some $C>0$ we have (analogously to the arguments used in \cite{Ueda})
\begin{eqnarray*}
\int   |  w(x,m) |^2 \, dx 
& =  &
\int   | \widehat w(\xi ,m) |^2 \, d\xi \geq  \int_{I}   | \widehat w(\xi ,m) |^2 \, d\xi \\
& \geq &
 \int_{I}   |  W(m,|\xi|^2 )\widehat \varphi (\xi )   |^2 \, d\xi \geq  \min_{\lambda \in [a,b]} | W(m,\lambda )|  \int_{I}   | \widehat \varphi (\xi )   |^2 \, d\xi \\
& \geq &
 C_\varphi   e^{ 2m\ln \min_{\lambda \in I } \mu_0(\lambda ) }  \qquad \forall m \in{\mathbb N},
\end{eqnarray*} 
and similarly for $ v(x,m) $.  
Now, if we take into account the cone of dependence, then for $q\geq 2$ we obtain
\begin{eqnarray*}
&  &
\|  w(x,m) \|_{L^q({\mathbb R}^n)}  \geq  C_\varphi (1+m)^{-\frac{q-2}{2q}n}\|  w(x,m) \|_{L^2({\mathbb R}^n)}  \geq 
  C_\varphi   e^{ \delta_\varphi  m}  \qquad \forall m \in{\mathbb N},\\ 
&  &
\|  v(x,m) \|_{L^q({\mathbb R}^n)}  \geq  C_\varphi (1+m)^{-\frac{q-2}{2q}n}\|  v(x,m) \|_{L^2({\mathbb R}^n)}  \geq 
  C_\varphi   e^{ \delta_\varphi  m}  \qquad \forall m \in{\mathbb N}.
\end{eqnarray*}
Theorem is proved. \hfill $\square$
\begin{corollary} 
\label{C4.4}
For all $m \in{\mathbb N} $ we have 
\begin{eqnarray*}
&  &
\max_{x \in {\mathbb R}^n }  |w(x,m)|  \geq   
  C_\varphi   e^{ \delta_\varphi  m}  ,\qquad 
\max_{x \in {\mathbb R}^n }  | v(x,m) | \geq   
  C_\varphi   e^{ \delta_\varphi  m}.
\end{eqnarray*}
\end{corollary}

\begin{theorem} 
\label{T4.5} 
Let  \,  $R(t)$\, be defined on \, ${\mathbb R}$\, non-constant,  positive, 
smooth function which is  $1$-periodic and $R'(0)=0 $.  
Then  for every function $ \varphi \in C_0^\infty ({\mathbb R}^n)$, which is different from the identical zero,
 there are positive numbers  $C_\varphi $, $\delta_\varphi  $   such that for $q \in [2,\infty]$ the solutions $w(x,t) $ and $v(x,t) $ 
of the equation 
\begin{eqnarray}
\label{LinEQ_v}
&  &
\partial_t^2  v + nR^{-1} R' \partial_t   v- R^{-2} \Delta _\gamma   v   =0
\end{eqnarray}
with the initial conditions 
\begin{eqnarray*} 
&  &
w(x,0)= \varphi (x )  ,\qquad w_t(x,0)=0,\\
\label{24_v}
&  &
 v(x,0)= 0,\qquad v_t(x,0)=\varphi (x  )   
\end{eqnarray*}
satisfy
\begin{eqnarray*}
&  &
\|  w(x,m) \|_{L^q({\mathbb R}^n)}  \geq  
  C_\varphi   e^{ \delta_\varphi  m}  \qquad \forall m \in{\mathbb N}, \\
&  &
\| v(x,m) \|_{L^q({\mathbb R}^n)} \geq
 C_\varphi   e^{ \delta_\varphi  m}   \qquad \forall m \in{\mathbb N}. 
\end{eqnarray*}
The constant  $\delta_\varphi  $ depends on the diameter of supp\,$\varphi  $ only. 
In particular, for all $m \in{\mathbb N} $ we have 
\begin{eqnarray*}
&  &
\max_{x \in {\mathbb R}^n }  |w(x,m)|  \geq   
  C_\varphi   e^{ \delta_\varphi  m}  ,\qquad 
\max_{x \in {\mathbb R}^n }  | v(x,m) | \geq   
  C_\varphi   e^{ \delta_\varphi  m}.
\end{eqnarray*}
\end{theorem}
\medskip

\section{The Completion of the Proof of the Main Theorem} 
\label{S6}

First we consider the case of $0<l<2$. We construct the wave map with the finite life-span via a solution of the wave equation and the geodesic,  which is the vertical line.
Let $\wt u^2 =  const >0$ be a second component of the constant (stationary)     wave map $(\wt u^1 ,\wt u^2 ) $. 
Consider the perturbation of its second initial data $\wt u^2_t (x,0)= 0  $, and look for the global  wave map  $(u^1 (x,t), u^2 (x,t))$, such that
\[
u^2 (x,0)= \wt u^2 ,\quad  u^2_t (x,0)= \varphi _1(x )\alpha \left( \frac{\alpha }{\mu }\varphi_0 (x) +\beta \right)^{\mu -1} ,\quad \varphi_0 , \varphi _1 \in C_0^\infty ({\mathbb R}^n)\,,
\]
$ \mu = 2/(2-l)$.
We set $\varphi_0 (x)=0 $ and $\beta  = \left[\wt u^2 \right]^{1/\mu } >0$.  
The integral (\ref{epsilonG2}) of the perturbed initial data
\begin{eqnarray*} 
&  &
\int_{\{0\}\times  {\mathbb R}^2}\sum_{   \gamma_0+ \gamma_1+\gamma_2+ \ldots  +\gamma_n =1,\ldots,G
\atop{\gamma _0=0,1}}  \left( \left[ \frac{\partial^\gamma  u^1}{(u^2)^{l/2}} \right]^2+  \left[ \frac{\partial^\gamma u^2}{(u^2)^{l/2}} \right]^2\right) dx \\
& = & 
\int_{\{0\}\times  {\mathbb R}^2} \sum_{ \gamma_1+\gamma_2+ \ldots  +\gamma_n \leq G -1 
\atop{\gamma_0 =1} } \left(   \frac{\partial^{\gamma} u^2}{(u^2)^{l/2}} \right)^2 dx = 
\int_{  {\mathbb R}^2} \sum_{  \gamma_1+\gamma_2+ \ldots  +\gamma_n \leq G -1}   \left(   \frac{ (\partial^{\gamma } \varphi _1(x )) \alpha ( \beta  )^{\mu -1} } {(\wt u^2)^{l/2} } \right)^2 dx 
\end{eqnarray*}
can be chosen arbitrarily small by appropriate choice of $\varphi _1 \in C_0^\infty ({\mathbb R}^n) $.

Then, by the uniqueness,  $  u^1 (x,t)=   \wt u^1  $,   and  for any real $\alpha \not=0 $ the function $v(x,t)$ defined by
\[
  v(x,t)  =  \frac{\mu }{\alpha }  \left( \left[u^2  (x,t)\right]^{1/\mu } - \beta \right)    , \qquad  u^2  (x,t)>0 , 
\]
is a solution of the linear equation  (\ref{3.16}) and takes initial values 
\begin{eqnarray*}
&  &
  v(x,0)  =0    , \qquad  v_t(x,0)  =      \varphi _1(x )  \,.
\end{eqnarray*}
Hence, for all $t \geq 0$ and all $x \in {\mathbb R}^n$ we have
\[
u^2 (x,t)=   \left( \frac{\alpha }{\mu } v(x,t) +\beta \right)^{\mu  } \,.
\]
According to Theorem~\ref{T4.5} with $q=\infty$ and initial condition (\ref{24}), for the appropriate choice of $\alpha  $ there exists a point $(t_{bp},x_{bp}) \in {\mathbb R}^{1+n} $  such that  
\begin{eqnarray*}
&  &
 \frac{\alpha }{\mu } v(t_{bp},x_{bp}) +\beta   \leq 0\,.
\end{eqnarray*} 
Thus, the solution blows up    in finite time.

In the case of $ l=2$ and $n=2$ we choose
\[
u^2 (x,0)= \wt u^2 ,\quad  u^2_t (x,0)= \wt u^2 \varphi _1(x )   ,\quad  \varphi _1 \in C_0^\infty ({\mathbb R}^n)\,.
\]
These new initial data have small integral (\ref{epsilonG2}) provided that the function  $\varphi _1 \in C_0^\infty ({\mathbb R}^n) $ has sufficiently small $L^2 ({\mathbb R}^n)$-norm.
Next we solve equation (\ref{LinEQ_v}) with the initial conditions
\begin{eqnarray*}
&  &
  v(x,0)=\ln \wt u^2, \quad  v_t(x,0)=\varphi _1(x )\,.
\end{eqnarray*}
Due to the uniqueness, the  vector valued function $( \wt u^1, e^{v(x,t)}) $ is a wave map. 
Theorem~\ref{T4.5} with $q=\infty$,  
for the appropriately chosen function $\varphi _1 $  completes the proof of the main theorem.
\hfill $\square$

\begin{remark}
The perturbation function $\varphi _1 \in C_0^\infty ({\mathbb R}^n) $ can be chosen spherically symmetric, arbitrarily small in every given  space
${\mathcal F}\supset  C_0^\infty ({\mathbb R}^n)$.
\end{remark}

\bigskip

\noindent
{\bf Acknowledgments}
\medskip

This paper was originated in November of 2009  when the second author was     visiting   Osaka University.  
He expresses  his gratitude   to the JSPS Grant-in-Aid for Scientific Research
(A) for   the financial support of the visit.

\end{document}